\documentclass[smallextended]{svjour3}

\smartqed

\newcommand{\be}{\begin{equation}}
\newcommand{\ee}{\end{equation}}
\newcommand{\beq}{\begin{eqnarray}}
\newcommand{\eeq}{\end{eqnarray}}
\newcommand{\nbeq}{\begin{eqnarray*}}
\newcommand{\neeq}{\end{eqnarray*}}

\journalname{Metrika}

\begin{document}

\title{Characterization of exponential distribution via regression of one
record value on two non-adjacent record values}
\author{George P. Yanev }
\institute{George P. Yanev \at
              Department of Mathematics, The University of Texas - Pan American\\
              1201 W. University Drive, Edinburg, Texas, 78539 USA\\
              Tel.: (956) 381-3632, Fax: (956) 381-2428, \email{yanevgp@utpa.edu}           }

\date{Received: date / Accepted: date}
\maketitle

\begin{abstract}
We characterize the exponential distribution as the only one which satisfies
a regression condition. This condition involves the regression function of a
fixed record value given two other record values, one of them being previous
and the other next to the fixed record value, and none of them are adjacent.
In particular, it turns out that the underlying distribution is exponential
if and only if given the first and last record values, the expected value of
the median in a sample of record values equals the sample midrange.

\end{abstract}



\titlerunning{Characterization via regression on two non-adjacent record
values} 


\institute{George P. Yanev \at
              Department of Mathematics, The University of Texas - Pan American\\
              1201 W. University Drive, Edinburg, Texas, 78539 USA\\
              Tel.: (956) 381-3632, Fax: (956) 381-2428, \email{yanevgp@utpa.edu}           }


\section{Introduction}

\label{intro} In 2006, on a seminar at the University of South Florida, Moe
Ahsanullah posed the question about characterizations of probability
distributions based on regression of a fixed record value with two
non-adjacent (at least two spacings away) record values as covariates. We
address this problem here.

To formulate and discuss  our results we need to introduce some notation as
follows. Let $X_1, X_2, \ldots$ be independent copies of a random variable $X
$ with absolutely continuous distribution function $F(x)$. An observation in
a discrete time series is called a (upper) record value if it exceeds all
previous observations, i.e., $X_j$ is a (upper) record value if $X_j>X_i$
for all $i<j$. If we define the sequence $\{T_n, n\ge 1\}$ of record times
by $T_1=1$ and $T_n=\min \{j:X_j>X_{T_{n-1}}, j>T_{n-1}\}$, $(n>1)$, then
the corresponding record values are $R_n=X_{T_n}$, $n=1,2,\ldots$ (see
Nevzorov (2001)).

Let $F(x)$ be the exponential distribution function
\begin{equation}  \label{exp_type}
F(x)=1-e^{\displaystyle -c(x-l_F)}, \qquad (x\ge l_F>-\infty),
\end{equation}
where $c>0$ is an arbitrary constant.  Let us mention that (\ref{exp_type})
with $l_F>0$  appears, for example, in reliability studies where $l_F$
represents the guarantee time; that is,  failure cannot occur before $l_F$
units of time have elapsed (see Barlow and Proschan (1996), p.13).

We study characterizations of exponential distributions in terms of the
regression of one record value with two other record values as covariates,
i.e., for $1\le k\le n-1$ and $r\ge 1$ we examine the regression function
\[
E[\psi(R_n)|R_{n-k}=u, R_{n+r}=v ], \qquad (v>u\ge l_F),
\]
where $\psi$ is a function that satisfies certain regularity conditions. Let $\overline{f}_{u,v}$ denote the average value of an integrable function $f(x)$ over the interval from $x = u$ to $x = v$, i.e.,
\[
\overline{f}_{u,v}=\frac{1}{v-u}\int_u^vf(t)dt.
\]
Yanev et al. (2008) prove, under some assumptions on the function $g$, that if $F$ is
exponential then for $1\le k\le n-1$ and $r\ge 1$,
\begin{equation}  \label{exp_property}
\hspace{-0.5cm}E\left[\frac{\displaystyle g^{(k+r-1)}(R_n)}{\displaystyle %
k+r-1}{\Big |}R_{n-k}=u, R_{n+r}=v\right] = {k-1+r-1 \choose k-1}\frac{\partial^{k+r-2}}{\partial u^{r-1}\partial v^{k-1}}\ \left(\overline{g^{\prime}}\!_{u,v}\right),  \nonumber
\end{equation}
where $v>u\ge l_F$ and $g^\prime$ is the derivative of $g$. Bairamov et al. (2005) study the particular case of (\ref%
{exp_property}) when both covariates are adjacent (one spacing away) to $R_n$%
. They prove, under some regularity conditions, that if $k=r=1$, then (\ref%
{exp_property}) is also sufficient for $F$ to be exponential. That is, $F$
is exponential if and only if
\[
E\left[g^{\prime }(R_n){\Big |}R_{n-1}=u, R_{n+1}=v \right]=\overline{g^{\prime}}\!_{u,v}, \qquad (v>u\ge l_F).
\]
Yanev et al. (2008) consider the case when only one of the two covariates is
adjacent to $R_n$ and show that, under some regularity assumptions, $F$ is
exponential if and only if (\ref{exp_property}) holds for $2\le k\le n-1$
and $r=1$, i.e.,
\[
E\left[\frac{g^{(k)}(R_n)}{k}{\Big |}R_{n-k}=u, R_{n+1}=v\right]= \frac{%
\partial^{k-1}}{\partial v^{k-1}}\ \left(\overline{g^{\prime}}\!_{u,v}\right), \qquad
(v>u\ge l_F).
\]
Here we address the case when both covariates are non-adjacent to $R_n$,
which turns to be more complex. Denote for $x\ge l_F$,
\[
H(x)=-\ln(1-F(x)) \qquad \mbox{and}\qquad h(x)=H^{\prime }(x),
\]
i.e., $H(x)$ is the cumulative hazard function of $X$ and $h(x)$ is its
hazard (failure) rate function. In this paper, under some additional
assumptions on the hazard rate $h(x)$ and the function $g(x)$, we extend the
results in Bairamov et al. (2005) to the case when both covariates are
non-adjacent. Namely, we shall prove that for fixed $2\le k \le n-1$ and $%
r\ge 2$, equation (\ref{exp_property}) is a necessary and sufficient
condition for $F(x)$ to be exponential. Note that the characterization for
the non-adjacent case given in Theorem 1B of Yanev et al. (2008) involves,
in addition to (\ref{exp_property}), one more regression condition. We shall
show here that (\ref{exp_property}) alone characterizes the exponential
distribution. This result provides a natural generalization of the known
special cases mentioned above. As a consequence of our main result, we
obtain Corollary 1 below, which seems to be of independent interest with
respect to possible statistical applications. Let us also mention that the
technique of our proof is different from that used by Dembi\'{n}ska and Weso{%
\l }owski (2000) in deriving characterization results in terms of regression
of a record value on another non-adjacent one.

Further on, for a given continuous  function $g(x)$ and positive integers $i$
and $j$, we denote
\begin{equation}  \label{M_not}
M(u,v)=\overline{g^{\prime}}\!_{u,v}=\frac{g(v)-g(u)}{v-u}, \ _iM_j(u,v)=\frac{\partial^{i+j}}{\partial
u^i\partial v^j}\ \left(M(u,v)\right), \quad (u\ne v),
\end{equation}
as well as $_iM(u,v)$ and $M_j(u,v)$ for the $i$th and $j$th partial
derivative of $M(u,v)$ with respect to $u$ and $v$, respectively.

\textbf{Theorem }\ \textit{Let $n$, $k$, and $r$ be integers, such that $%
2\le k\le n-1$ and $r\ge 2$. Assume that $F(x)$ satisfies the following
conditions. }

\textit{(i) The $n$th derivative $F^{ (n)}(x)$ where $n=\max\{k,r\}$ is
continuous in $(l_F,\infty)$; }

\textit{(ii) $h(x)$ is nowhere constant in a small interval $(l_F,
l_F+\varepsilon)$ for $\varepsilon>0$; }

\textit{(iii) $h(l_F+)>0$ and $\left|h^{(n)}(l_F+)\right|<\infty$ for $n\le \max(2,r-1)$. }

\textit{\noindent Suppose the function $g(x)$ satisfies }

\textit{(iv) $g(x)$ is continuous in $(l_F, \infty)$ and $g^{(k+r-1)}(x)$ is
continuous in $(l_F,\infty)$; }

\textit{(v) $_{r-1}M_k(l_F+,v)\ne 0$ for $v\ge l_F$; }

\textit{(vi) if $r=2$ then
$|g^{(k+2)}(l_F+)|<\infty$, and if $r\ge 3$ then $|g^{(k+2r-1)}(l_F+)|<\infty$. }

\textit{\noindent Then (\ref{exp_property}) holds if and only if $X$ has the
exponential distribution (\ref{exp_type}) with $c=h(l_F+)$.  }

\textbf{Remark}. I conjecture that the assumption (vi) can be weakened to $%
|g^{(k+r)}(l_F+)|<\infty$ for any $r\ge 2$, retaining the symmetry with
respect to $k$ and $r$ from the case $r=2$. One can verify this in the case $%
r=3$ by extending the approximation formula in Lemma~4.

We refer to Leemis (1995) for distributions, related to reliability and
lifetime modeling, whose hazard functions satisfy the assumptions (ii) and
(iii). Also the two corollaries below provide examples of functions $g(x)$
which satisfy the assumptions of the Theorem.

We continue with two interesting particular choices for 
$g(x)$. First, setting
\[
g(x)=\frac{x^{k+r}}{(k+r)!} \quad \mbox{and thus}\quad \frac{g^{(k+r-1)}(x)}{%
k+r-1}=\frac{x}{k+r-1},
\]
one can see that the assumptions (iv)-(vi) of the Theorem are satisfied and
\[
{k+r-2 \choose k-1}\ _{r-1}M_{k-1}(u,v) =\frac{1}{k+r-1}\frac{ru+kv}{k+r}.
\]
Therefore, we obtain the following corollary.

\textbf{Corollary 1}\ \textit{Let $n$, $k$, and $r$ be integers, such that $%
2\le k\le n-1$ and $r\ge 2$. Suppose assumptions (i)-(iii) of the Theorem
hold. Then $X$ has the exponential distribution (\ref{exp_type}) with $%
c=h(l_F+)$ 
if and only if
\begin{equation}  \label{thmg}
E[R_n|R_{n-k}=u, R_{n+r}=v]=\frac{ru+kv}{k+r}, \qquad (v>u\ge l_F).
\end{equation}
}

Note that the right-hand side of (\ref{thmg}) is a weighted average of the
two covariate values - each covariate being given weight proportional to the
number of spacings $R_n$ is away from the other covariate.
In particular, (\ref{thmg}) with $k=r$ becomes 
\[
E[R_n|R_{n-k}=u, R_{n+k}=v]= \frac{u+v}{2},\qquad (2\le k \le n-1).
\]
This last equation allows the following interpretation. Suppose we observe $%
2n-1$ record values $R_1, \ldots , R_{2n-1}$ where $n\ge 2$. Then $X$ is
exponential if and only if, given the first and last record values, the
expected value of the median $R_n$ in the sample equals the sample midrange.

We continue with another choice of $g(x)$ from (\ref{exp_property}). Let $%
l_F>0$ and
\[
g(x)=\frac{(-1)^{k+r-1}}{(k+r-1)!}\ \frac{1}{x}\quad \mbox{and thus}\quad
\frac{g^{(k+r-1)}(x)}{k+r-1}=\frac{1}{(k+r-1)x^{k+r}}.
\]
It is not difficult to see that the assumptions (iv)-(vi) of the Theorem are
satisfied and
\[
{k+r-2 \choose k-1}\ _{r-1}M_{k-1}(u,v) =\frac{1}{(k+r-1)u^rv^k}.
\]
Hence, the Theorem implies the following result.

\textbf{Corollary 2}\ \textit{Let $n$, $k$, and $r$ be integers, such that $%
2\le k\le n-1$ and $r\ge 2$. Suppose assumptions (i)-(iii) of the Theorem
hold.
Then $X$ has the exponential distribution (\ref{exp_type}) with $c=h(l_F+)$
if and only if
\begin{equation}  \label{thmg3}
E\left[\frac{1}{R^{k+r}_n}{\Big |}R_{n-k}=u, R_{n+r}=v\right]=\frac{1}{u^rv^k%
}, \ \ (v>u\ge l_F> 0).
\end{equation}
}

Finally, let us mention that, following Bairamov et al. (2005), one can
obtain an extension of the Theorem that involves monotone transformations of
$X$, see also Yanev et al. (2008), Theorem 3. Consequently, the
characterization examples given in the above two papers  can be modified for
the case of non-adjacent covariates.

\section{Preliminaries}

\label{sec:1}

In this section we present four technical lemmas, which we use in Section 3
to prove the Theorem.  First, we prove an identity that links the
derivatives of $g(x)$ with those of $M(u,v)=(g(v)-g(u))/(v-u)$.  Denote $%
(n)_{(m)}=n(n-1)\ldots (n-m+1)$ $(m\ge 1)$; $n_{(0)}=1$.

\vspace{0.3cm}\textbf{Lemma 1}\ For any positive integer $k$ and $n\ge 2$

\noindent
\begin{equation}  \label{le1}
(n-1)!g^{(k+n-1)}(v)=\sum_{i=0}^{n}{n \choose i}(k+n-1)_{(n-i)}(v-u)^i\
_{n-1}M_{k-1+i}(u,v), \qquad (v>u).
\end{equation}
\textbf{Proof}.\ For simplicity write $_iM_j$ for $_iM_j(u,v)$. According to
Lemma~1 in Yanev et al. (2008), we have for $i,j\ge 1$

\begin{equation}  \label{lemma_old}
g^{(j)}(v)=(v-u)M_j+jM_{j-1}, \quad i\ _{i-1}M_j=(v-u)\ _iM_j+j\ _iM_{j-1},
\qquad (v>u).
\end{equation}
To prove (\ref{le1}) we use induction with respect to $n$. Referring to (\ref%
{lemma_old}), we have
\begin{eqnarray}  \label{lemma_1_i}
g^{(k+1)}(v) & = & (v-u)M_{k+1}+(k+1)M_k  \nonumber \\
& = & (v-u)[(v-u)\ _1M_{k+1}+(k+1)\ _1M_k]+(k+1)[(v-u)\ _1M_k+k\ _1M_{k-1}]
\nonumber \\
& = & (k+1)k\ _1M_{k-1} +2(k+1)(v-u)\ _1M_k+(v-u)^2\ _1M_{k+1},  \nonumber
\end{eqnarray}
which is (\ref{le1}) with $n=2$. To complete the proof,  assuming (\ref{le1}%
), we need to show that
\begin{equation}  \label{le1_n+1}
n!g^{(k+n)}(v)=\sum_{i=0}^{n+1}{n+1 \choose i}(k+n)_{(n+1-i)}(v-u)^i\
_nM_{k-1+i}.
\end{equation}
Differentiating both sides of (\ref{le1}) with respect to $v$ and
multiplying by $n$, we obtain
\begin{equation}  \label{diff_formula}
n!g^{(k+n)}(v)= \sum_{i=0}^{n}{n \choose i}(k+n-1)_{(n-i)}n\left[
i(v-u)^{i-1}\ _{n-1}M_{k+i-1}+(v-u)^i\ _{n-1}M_{k+i}\right].
\end{equation}
Applying the second formula in (\ref{lemma_old}) repeatedly, we have
\begin{eqnarray}  \label{diff_form_2}
\lefteqn{n\left[ i(v-u)^{i-1}\ _{n-1}M_{k+i-1}+(v-u)^i\ _{n-1}M_{k+i}\right]}
\nonumber \\
& = & i(v-u)^{i-1}\left[(v-u)\ _{n}M_{k+i-1}+(k+i-1)\ _{n}M_{k+i-2}\right]
\nonumber \\
& & +(v-u)^i\left[(v-u)\ _{n}M_{k+i}+ (k+i)\ _{n}M_{k+i-1}\right] \\
& = & (v-u)^{i+1}\ _{n}M_{k+i}+(k+2i)(v-u)^i\
_{n}M_{k+i-1}+i(k+i-1)(v-u)^{i-1}\ _{n}M_{k+i-2}.  \nonumber
\end{eqnarray}
Therefore, by (\ref{diff_formula}) and (\ref{diff_form_2}), we have
\begin{eqnarray}  \label{new_sum}
n!g^{(k+n)}(v) & = & \sum_{i=0}^{n}{n \choose i}(k+n-1)_{(n-i)}(v-u)^{i+1}%
\ _{n}M_{k+i}  \nonumber \\
& & + \sum_{i=0}^{n}{n \choose i}(k+n-1)_{(n-i)}(k+2i)(v-u)^i\
_{n}M_{k+i-1}  \nonumber \\
& & +\sum_{i=0}^{n}{n \choose i}(k+n-1)_{(n-i)}i(k+i-1)(v-u)^{i-1}\
_{n}M_{k+i-2} \\
& = & S_1 + S_2 + S_3, \qquad \mbox{say}.  \nonumber
\end{eqnarray}
Changing the summation index to $l=i+1$ we obtain
\begin{equation}  \label{sum_1}
S_1 = \sum_{l=0}^{n+1}{n \choose l-1}(k+n-1)_{(n-l+1)}(v-u)^l\
_{n}M_{k+l-1}  \nonumber
\end{equation}
and setting $l=i-1$, we have
\begin{equation}  \label{sum_2}
S_3= \sum_{l=0}^{n-1} { n \choose l+1} (l+1)(k+n-1)_{(n-l-1)}(k+l)(v-u)^l\
_{n}M_{k+l-1},  \nonumber
\end{equation}
assuming ${n \choose l}=0$ for $l=-1$ or $l>n$. Now, observing that
\begin{eqnarray}  \label{identity}
\lefteqn{\hspace{-1.5cm}{n \choose i-1}(k+n-1)_{(n-i+1)}+{ n \choose
i}(k+n-1)_{(n-i)}(k+2i)+{n \choose i+1}(k+n-1)_{(n-i-1)}(i+1)(k+i)}
\nonumber \\
& = & (k+n-1)_{(n-i)}\left[{n \choose i-1}(k+i-1)+{n \choose i}(k+2i)+{n \choose i+1}(i+1)\right] \\
& = & {n+1 \choose i}(k+n)_{(n-i+1)},  \nonumber
\end{eqnarray}
one can see that (\ref{new_sum})-(\ref{identity}) imply (\ref{le1_n+1})
which completes the proof of the lemma. 

For simplicity, further on we denote, for integer $i,j\ge 0$ and $v\ge l_F$,
\[
\ _iM_j(v)=\! \ _iM_j(l_F+,v).
\]
The following result holds.

\vspace{0.3cm}\textbf{Lemma 2}\ \textit{If $|g^{(i+j+1)}(l_F+)|<\infty$ for
any non-negative integers $i$ and $j$, then
\begin{equation}  \label{M_at l_F}
\lim_{v\to l_F+}{i+j \choose i} \ _iM_j(v)=\frac{g^{(i+j+1)}(l_F+)}{i+j+1}.
\end{equation}
Also, if $|g^{(i+j+1-m)}(l_F+)|<\infty$ for $m=1,2,\ldots$, then
\begin{equation}  \label{M_at_l_F_extra}
\lim_{v\to l_F+}(v-l_F)^m\ _iM_j(v)=0
\end{equation}
}

\textbf{Remark}. Note that for $i=k-1$ and $j=r-1$, (\ref{M_at l_F}) implies
that the limit of the right-hand side of (\ref{exp_property}) as $v\to l_F+$
equals $g^{(k+r-1)}(l_F+)/(k+r-1)$.

\textbf{Proof}.\ We use induction with respect to the sum $i+j$. Clearly $%
\lim_{v\to \l_F+}M(v)=g^{\prime }(l_F+)$. Applying L'Hopital's rule, we
have
\begin{eqnarray*}
\lim_{v\to l_F+}M_1(v) & = & \lim_{v\to l_F+}\frac{g^{\prime }(v)-M(v)}{v-l_F%
} \\
& = & \lim_{v\to l_F+}g^{\prime \prime }(v)-\lim_{v\to l_F+}M_1(v).
\end{eqnarray*}
Hence, $\lim_{v\to l_F+}M_1(v)=g^{\prime \prime }(l_F+)/2$. Similarly, $%
\lim_{v\to l_F+}\ _1M(v)=g^{\prime \prime }(l_F+)/2$. This verifies (\ref%
{M_at l_F}) for $i+j=0$ and $i+j=1$. Assuming that (\ref{M_at l_F}) is true
for $0\le i+j\le n$, we will prove it for $i+j=n+1$. By the second equation
in (\ref{lemma_old}) and L'Hopital's rule (the numerator below
approaches zero by the induction assumption) we have
\begin{eqnarray*}
\lim_{v\to l_F+}\ _iM_j(v) & = & \lim_{v\to l_F+} \frac{i\ _{i-1}M_j(v)-j\
_iM_{j-1}(v)}{v-l_F} \\
& = & \lim_{v\to l_F+} i\ _{i-1}M_{j+1}(v)-j \lim_{v\to l_F+}\ _iM_j(v).
\end{eqnarray*}
That is,
\[
\lim_{v\to l_F+}\ _iM_j(v)=\frac{i}{j+1}\lim_{v\to l_F+}\ _{i-1}M_{j+1}(v).
\]
Iterating, we obtain
\begin{equation}  \label{one_sided}
\lim_{v\to l_F+}\ _iM_j(v)= \frac{i!j!}{(i+j)!}\lim_{v\to l_F+}\ M_{j+i}(v).
\end{equation}
Now, by the first equation in (\ref{lemma_old}) and L'Hopital's rule
(the numerator below approaches zero by the induction assumption) we have
\begin{eqnarray*}
\lim_{v\to l_F+}M_{i+j}(v) & = & \lim_{v\to l_F+} \frac{%
g^{(i+j)}(v)-(i+j)M_{i+j-1}(v)}{v-l_F} \\
& = & \lim_{v\to l_F+}g^{(i+j+1)}(v) - (i+j)\lim_{v\to l_F+}M_{i+j}(v)
\end{eqnarray*}
and hence
\[
\lim_{v\to l_F+}M_{i+j}(v)=\frac{1}{i+j+1}g^{(i+j+1)}(l_F+).
\]
Substituting this into (\ref{one_sided}) we complete the proof of the
induction step.

Let us now prove (\ref{M_at_l_F_extra}). Using induction and the second
equation in (\ref{lemma_old}), it is not difficult to see that for $%
m=0,1,\ldots$
\[
(v-l_F)^m \ _iM_j(v)=\sum_{k=0}^m {m \choose k}(-1)^k i_{(m-k)}j_{(k)}\
_{i-m+k}M_{j-k}(v).
\]
Passing to the limit as $v\to l_F+$ and applying (\ref{M_at l_F}) we find
\begin{eqnarray*}
\lim_{v\to l_F+} (v-l_F)^m\ _iM_j(v) & = & \sum_{k=0}^m {m \choose k}%
(-1)^k i_{(m-k)}j_{(k)}\lim_{v\to l_F+}\ _{i-m+k}M_{j-k}(v) \\
& = & \sum_{k=0}^m {m \choose k}(-1)^k i_{(m-k)}j_{(k)}\frac{(i-m+k)!(j-k)!%
}{(i+j+1-m)!}g^{(i+j+1-m)}(l_F+) \\
& = & \frac{i!j!}{(i+j+1-m)!}g^{(i+j+1-m)}(l_F+)\sum_{k=0}^m {m \choose k}%
(-1)^k \\
& = & 0.
\end{eqnarray*}
The proof of the lemma is complete.

The next lemma establishes some identities and limit results involving
\begin{equation}  \label{new_var}
w(v)=\frac{h(v)}{H(v)}(v-l_F),\qquad (v>l_F).
\end{equation}

\vspace{0.3cm}\textbf{Lemma 3}\ \textit{For $v> l_F$,
\begin{eqnarray}  \label{frac_h}
\hspace{-0.5cm}\frac{h^{\prime }(v)}{h(v)}(v-l_F) & = & \frac{w^{\prime }(v)%
}{w(v)}(v-l_F)+w(v)-1.
\end{eqnarray}
If $F^{\prime \prime }(v)$ is continuous in $(l_F,\infty)$, $h(l_F+)>0$, and
$h^{\prime }(l_F+)\ne 0$, then
\begin{equation}  \label{fraction_lim_1}
\lim_{v\to l_F+}w(v)=1,\quad \lim_{v\to l_F+}\frac{(v-l_F)w^{\prime }(v)}{%
w(v)-1}=1,
\end{equation}
and
\begin{equation}  \label{m_frac}
\lim_{v\to l_F+}\frac{(v-l_F)^2}{w(v)-1}=0.
\end{equation}
}

\textbf{Proof}.\ Differentiating (\ref{new_var}) with respect to $v$, it is
not difficult to obtain (\ref{frac_h}).
Applying L'Hopital's rule, we obtain that as $v\to l_F+$
\begin{eqnarray}  \label{w_prime_2}
w^{\prime }(v) & = & \frac{[H(v)h^{\prime}(v)-h^2(v)](v-l_F)+H(v)h(v)}{H(v)^2} \\
& \sim & \frac{[H(v)h^{\prime \prime }(v)-h(v)h^{\prime
}(v)](v-l_F)+2H(v)h^{\prime }(v)}{2H(v)h(v)}  \nonumber \\
& \to & \frac{h^{\prime }(l_F+)}{2h(l_F+)}\ne 0 .  \nonumber
\end{eqnarray}
Now, the continuity of $w(v)$ implies that $w(l_F+)=1$. It follows by the
mean-value theorem and (\ref{w_prime_2}) that
\[
\lim_{v\to l_F+}\frac{(v-l_F)w^{\prime }(v)}{w(v)-1} = \lim_{v\to l_F+}\frac{%
w^{\prime }(v)}{w^{\prime }(\eta)}=1, \qquad (l_F<\eta<v),
\]
i.e., the second limiting result in (\ref{fraction_lim_1}). Finally,
applying L'Hopital's rule, it is not difficult to obtain (\ref{m_frac}).
The proof of the lemma is complete.

For positive integers $n$, $r$, and $k$, define the sequence $%
\{d_n(v)\}_{n=1}^\infty$ for $v>l_F$ by the recurrence
\begin{equation}  \label{d_seq}
d_1(v)=\frac{d}{d v}\{\ _{r-1}M_{k-1}(v)H^{n+k-1}(v)\} \quad \mbox{and}\quad
d_{n+1}(v)=\frac{d}{d v}\left\{\frac{d_n(v)}{h(v)}\right\}.
\end{equation}
In the lemma below, we derive an expansion of $d_n(v)$ in terms of $\
_{r-1}M_j(v)$ and $H^j(v)$ for $k-1\le j\le n+k-1$. Note that if $k>i$, then
${i\choose k}=0$ and $\sum_{j=k}^i(\cdot )=0$.

\vspace{0.3cm}\textbf{Lemma 4}\ \textit{The following identity is true for $%
n=1,2,\ldots$
\begin{eqnarray}  \label{d_formula}
d_n(v)& = & \sum_{j=0}^2 {n \choose j} (k+n-1)_{(n-j)}\ _{r-1}M_{k-1+j}(v)%
\frac{H^{k-1+j}(v)}{h^{j-1}(v)}  \nonumber \\
& & -{n \choose 2}(k+n-1)_{(n-2)}\ _{r-1}M_k(v) \frac{h^{\prime k+1}(v)}{%
h^2(v)}+ \sum_{j=3}^{n} c_j(v)H^{k-1+j}(v),
\end{eqnarray}
provided that the left and right-hand sides are well-defined. }

\textit{If $h(l_F+)\ne 0$, $|h^{(n-1)}(l_F+)|<\infty$ and $%
|g^{(k+r+n-1)}(l_F+)|<\infty$ for $n=3,4,\ldots$, then
\begin{equation}  \label{d_formula2}
\limsup_{v\to l_F+}\left|\sum_{j=3}^{n} c_j(v)H^{k-1+j}(v)\right|<\infty.
\end{equation}
}

\textbf{Proof}.\ Using induction, one can prove that for $n=1,2,\ldots$
\begin{eqnarray}  \label{ind}
d_n(v) & = & -\frac{h^{\prime }(v)}{h^2(v)}d_{n-1}(v)+\frac{1}{h(v)}%
d^{\prime }_{n-1}(v) \\
& = & \sum_{j=1}^n c_{j,n}(v)\frac{d^{n-j}}{dv^{n-j}} d_1(v)  \nonumber
\end{eqnarray}
and $c_{j,n}(v)$ satisfy the following equations for $j=2,3,\ldots, n$,
\[
c_{j,n}(v) = \frac{1}{h(v)}c^{\prime }_{j-1,n-1}(v)-\frac{h^{\prime }(v)}{%
h^2(v)}c_{j-1,n-1}(v)+\frac{1}{h(v)}c_{j,n-1}(v),
\]
where $c_{j,i}(v)=0$ if $j>i$ and $c_{1,n}(v)=1/h^{n-1}(v)$. It is not
difficult to obtain
\[
c_{1,n}(v)=\frac{1}{h^{n-1}(v)}, \qquad c_{2,n}(v)=-{n \choose 2}\frac{%
h^{\prime }(v)}{h^n(v)},
\]
and
\[
c_{3,n}(v)={n \choose 3}\left[ \frac{3(n+1)}{4}\frac{(h^{\prime}(v))^2}{%
h^{n+1}(v)}-\frac{h^{\prime \prime }(v)}{h^n(v)}\right].
\]
Note that 
$|c_{j,n}(v)|<\infty$ if $h(v)\ne~0$ and $|h^{(j-1)}(v)|<\infty$ for $1\le
j\le n$ .

For simplicity, further on in the proof we drop the left subscript $r-1$ in $%
\ _{r-1}M_j(v)$ and write $M_j(v)$ instead. Using Leibniz rule for
differentiation of the product of two functions, we have for $m\ge 1$
\begin{eqnarray}  \label{product_rule}
\frac{d^{m-1}}{d v^{m-1}}d_1(v) & = & \frac{d^{m}}{d v^{m}}%
\left\{M_{k-1}(v)H^{n+k-1}(v)\right\}  \nonumber \\
& = & \sum_{j=0}^{m}{m \choose j}M_{k-1+j}(v)\frac{d^{m-j}}{d v^{m-j}}%
H^{n+k-1}(v)  \nonumber
\end{eqnarray}
and hence (\ref{ind}) becomes
\begin{eqnarray}  \label{ind3}
d_n(v) & = & c_{1,n}(v)  \nonumber \\
& & \hspace{-1.2cm} \times \left[M_{k-1}(v)\frac{d^{n}}{d v^{n}}%
H^{n+k-1}(v)+nM_k(v)\frac{d^{n-1}}{d v^{n-1}}H^{n+k-1}(v)+{n \choose 2}%
M_{k+1}(v)\frac{d^{n-2}}{d v^{n-2}}H^{n+k-1}(v)\right]  \nonumber \\
& & \hspace{-1cm}+c_{2,n}(v)\left[M_{k-1}(v)\frac{d^{n-1}}{d v^{n-1}}%
H^{n+k-1}(v)+(n-1)M_k(v)\frac{d^{n-2}}{d v^{n-2}}H^{n+k-1}(v)\right]
\nonumber \\
& & \hspace{-1cm}+c_{3,n}(v)M_{k-1}(v)\frac{d^{n-2}}{d v^{n-2}}H^{n+k-1}(v)
+ S(v,M,H), \qquad \mbox{say}.
\end{eqnarray}
The last term, $S(v,M,H)$, in (\ref{ind}) does not include derivatives of $%
H^{n+k-1}(v)$ of order higher than $n-3$ and it is given by
\begin{eqnarray*}
S(v,M,H) & = & \sum_{j=0}^2\left[\sum_{i=j+1}^{n-2+j}{n-2+j \choose i}%
M_{k-1+i}(v)\frac{d^{n-2+j-i}}{dv^{n-2+j-i}}H^{n+k-1}(v)\right]c_{3-j,n}(v)
\\
& & + \sum_{j=3}^n\left[ \sum_{i=0}^{n-j}{n-j \choose i}M_{k-1+i}(v)\frac{%
d^{n-j-i}}{dv^{n-j-i}}H^{n+k-1}(v)\right]c_{j+1,n}(v).
\end{eqnarray*}
Note that $|S(v,M,H)|<\infty$ if $|M_{k-1+n}(v)|<\infty$ and $%
|c_{j,n}(v)|<\infty$ for $0\le j\le n$.

Recall the formula for the $n$th derivative of $f^m(v)$ for positive integer
$m$ (e.g., Wolfram Research~(2009)).
\[
\frac{d^{n}}{d v^{n}}f^m(v) = \sum_{i_1=0}^n\sum_{i_2=0}^{n-i_1}\ldots
\sum_{i_{m-1}=0}^{n-\sum_{j=1}^{m-2}i_j} \left(\prod_{p=1}^{m-1}{
n-\sum_{j=1}^{p-1}i_j \choose i_p}\right) \left(\prod_{j=1}^{m}\frac{d^{{i_j}}}{%
d v^{i_j}}f(v)\right),%
\]
where $i_1,i_2,\ldots,i_m$ is a partition of $n$. Observe that a term in
the right-hand side includes $f^j(v)$ if exactly $j$ of $i_1, \ldots, i_m$
are zeros. Let us apply this formula to $f^m(v)=H^{n+k-1}(v)$. Setting $%
m=n+k-1$, we see that there are at least $k-1$ zeros in the partition $%
i_1,i_2,\ldots,i_{n+k-1}$. Also the positions of $j$ zeros among the terms
of the partition $i_1,i_2,\ldots,i_{n+k-1}$ can be selected in ${n+k-1 \choose j}$ ways.
Therefore, we can list the terms in the right-hand side,
starting with the one that contains $H^{k-1}(v)$, as follows.
\begin{eqnarray*}  \label{gen_der}
\frac{d^{n}}{d v^{n}}H^{n+k-1}(v) & = & {n+k-1 \choose k-1}{n \choose 1}%
\ldots {1  \choose 1}(H^{\prime}(v))^ nH^{k-1}(v) \\
& & \hspace{-2cm}+{n-1 \choose 1}{n+k-1 \choose k}{n \choose 2}{
n-2 \choose 1}{n-3 \choose 1}\ldots {1 \choose 1}H^{\prime \prime
}(v)(H^{\prime}(v))^{n-2}H^k(v) \\
& & \hspace{-2cm}+{ n-2 \choose 2}{n+k-1 \choose k+1}{n \choose 2}{%
n-2 \choose 2}{n-4 \choose 1}{n-5 \choose 1} \ldots {1 \choose 1}%
(H^{\prime \prime}(v))^{2}(H^{\prime}(v))^{n-4}H^{k+1}(v) \\
& & \hspace{-2cm} +{n-2 \choose 1}{n+k-1 \choose k+1}{n \choose 3}{n-3 \choose 1}{n-4 \choose 1}
\ldots {1 \choose 1}H^{\prime \prime
\prime }(v)(H^{\prime}(v))^{n-3}H^{k+1}(v) \\
& & \hspace{-2cm}+\sum_{j=3}^{n-1}c_{j}(v,n)H^{k-1+j}(v)  \nonumber \\
& & \hspace{-2.5cm} =(k+1)ka_kh^n(v)H^{k-1}(v)+{n \choose 2}%
(k+1)a_kh^{n-2}(v)h^{\prime k}(v) \\
& & \hspace{-2cm}+{n \choose 3}a_k\left[\frac{3(n+1)}{4}%
h^{n-4}(v)(h^{\prime 2}+h^{n-3}(v)h^{\prime \prime }(v)\right]H^{k+1}(v)
+\sum_{j=3}^{n-1}c_{j}(v,n)H^{k-1+j}(v),  \nonumber
\end{eqnarray*}
where $a_k=(n+k-1)!/(k+1)!$ and $c_{j}(v,n)$ are functions of $h(v)$ and its
derivatives. Note that $|c_{j}(v,n)|<\infty$ if $|h^{(j)}(v)|<\infty$ for $%
j=3,\ldots n-1$. Similarly, for the derivatives of $H^{n+k-1}(v)$ of order $%
n-1$ and $n-2$ we find
\begin{eqnarray*}
\frac{d^{n-1}}{d v^{n-1}}H^{n+k-1}(v) & = & (k+1)a_kh^{n-1}(v)H^k(v)+{
n-1 \choose 2}a_kh^{n-3}(v)h^{\prime k+1}(v) \\
& & +\sum_{j=3}^{n-1}c_{j}(v,n-1)H^{k-1+j}(v),
\end{eqnarray*}
where $|c_{j}(v,n-1)|<\infty$ if $|h^{(j-1)}(v)|<\infty$ for $j=3,\ldots n-1$%
; and
\begin{eqnarray*}
\hspace{-3.7cm}\frac{d^{n-2}}{d v^{n-2}}H^{n+k-1}(v)& = &
a_kh^{n-2}(v)H^{k+1}(v)+\sum_{j=3}^{n-1}c_{j}(v,n-2)H^{k-1+j}(v),
\end{eqnarray*}
where $|c_{j}(v,n-2)|<\infty$ if $|h^{(j-2)}(v)|<\infty$ for $j=3,\ldots n-1$%
. Using the above three formulas we write (\ref{ind3}) as
\begin{eqnarray*}
\lefteqn{d_n(v)} \\
& = & \frac{a_kM_{k-1}(v)}{h^{n-1}(v)}\left[ h^n(v)H^{k-1}(v)+{n \choose 2}%
(k+1)h^{n-2}(v)h^{\prime k}(v)\right. \\
& & +\left. 3{n \choose 4}h^{n-4}(v)(h^{\prime 2}H^{k+1}(v)+{n \choose 3}%
h^{n-3}(v)h^{\prime \prime k+1}(v)+\sum_{j=3}^{n-1}c_{j}(v,n)H^{k-1+j}(v)%
\right] \\
& & + \frac{na_kM_k(v)}{h^{n-1}(v)}\left[(k+1)h^{n-1}(v)H^k(v)+{n-1 \choose 2}h^{n-3}(v)h^{\prime k+1}(v)\right. \\
& & \left.+\sum_{j=3}^{n-1}c_{j}(v,n-1)H^{k-1+j}(v)\right] + \frac{a_kM_{k+1}%
}{h^{n-1}(v)}{n \choose 2} h^{n-2}(v)H^{k+1}(v)+%
\sum_{j=3}^{n-1}c_{j}(v,n-2)H^{k-1+j}(v) \\
& & -{n \choose 2}\frac{a_kM_{k-1}(v)h^{\prime }(v)}{h^n(v)}\left[
(k+1)h^{n-1}(v)H^k(v)+{n-1 \choose 2}h^{n-3}(v)h^{\prime k+1}(v)\right. \\
& & \left.+\sum_{j=3}^{n-1}c_{j}(v,n-1)H^{k-1+j}(v)\right] -{n \choose 2}%
\frac{a_k(n-1)M_{k}(v)h^{\prime }(v)}{h^n(v)}h^{n-2}(v)H^{k+1}(v) \\
& & +\sum_{j=3}^{n-1}c_{j}(v,n-2)H^{k-1+j}(v)+ {n \choose 3}\left[ \frac{%
3(n+1)}{4}\frac{(h^{\prime 2}}{h^{n+1}(v)}-\frac{h^{\prime \prime }(v)}{%
h^n(v)}\right]a_kM_{k-1}(v)h^{n-2}(v)H^{k+1}(v) \\
& & +\sum_{j=3}^{n-1}c_{j}(v,n-2)H^{k-1+j}(v) + S(v,M,H) \\
& = & k(k+1)a_kM_{k-1}(v)h(v)H^{k-1}(v)+n(k+1)a_kM_k(v)H^k(v)+{n \choose 2}%
\frac{a_k}{h(v)}M_{k+1}(v)H^{k+1}(v) \\
& & -{n \choose 2}\frac{a_kh^{\prime }(v)}{h^2(v)}M_k(v)H^{k+1}(v)+%
\sum_{j=3}^{n-1} b_j(v)H^{k-1+j}(v) + S(v,M,H),
\end{eqnarray*}
where $|b_j(v)|<\infty$ if $h(v)\ne 0$, $|h^{(j)}(v)|<\infty$, and $%
|M_{k+1}(v)|<\infty$. This is equivalent to (\ref{d_formula}). The statement
in (\ref{d_formula2}) follows from the conditions for finiteness of $%
\sum_{j=3}^{n-1} b_j(v)H^{k-1+j}(v)$ and $S(v,M,H)$ given in the proof above.


\section{Proof of the Theorem}

\label{sec:2} 
It follows from Lemma 2 in Yanev et al. (2008) that (\ref{exp_property}) is
a necessary condition for $X$ to be exponential. Here we shall prove the
sufficiency. The scheme of the proof is as follows: (i) differentiate (\ref%
{exp_property}) $r$ times with respect to $v$, to obtain a differential
equation for $H(v)$; (ii) make an appropriate change of variables; (iii)
assuming that there is a non-exponential  solution, reach a contradiction.

Recall the formula for the conditional density $f_{k,r}(t|u,v)$, say, of $R_n
$ given $R_{n-k}=u$  and $R_{n+r}=v$, where $1\leq k<n$ and $r\geq 1$.
Namely, it can be derived using the Markov property of record values (e.g.,
Ahsanullah (2004), p.6) that for $u<t<v$
\begin{equation}  \label{cond_den}
f_{k,r}(t|u,v) = \frac{\displaystyle (k+r-1)!}{\displaystyle (k-1)!(r-1)!}%
\frac{\displaystyle (H(t)-H(u))^{k-1}(H(v)-H(t))^{r-1}}{\displaystyle %
(H(v)-H(u))^{k+r-1}}H^{\prime }(t).
\end{equation}
Using (\ref{cond_den}) 
we can write (\ref{exp_property}) as
\[
\int_u^vg^{(k+r-1)}(t) (H(t)-H(u))^{k-1}(H(v)-H(t))^{r-1}dH(t)= \
_{r-1}M_{k-1}(u,v)(H(v)-H(u))^{k+r-1}.
\]
The continuity of $F(x)$ implies $H(l_F+)=0$ and hence, letting $u\to l_F+$,
we have
\[
\int_{l_F}^vg^{(k+r-1)}(t) H^{k-1}(t)(H(v)-H(t))^{r-1}dH(t)= \
_{r-1}M_{k-1}(v)H^{k+r-1}(v).
\]
Differentiating the above equation $r$ times with respect to $v$, dividing
by $h(v)>0$ prior to every differentiation (after the first one), and
applying Lemma 4 with $n=r\ge 2$, we obtain
\begin{eqnarray}  \label{first_eqn}
\lefteqn{\hspace{-1.5cm}(r-1)!g^{(k+r-1)}(v)h(v)H^{k-1}(v)=d_r(v)} \\
& = & \sum_{j=0}^2 {r \choose j} (k+r-1)_{(r-j)}\ _{r-1}M_{k-1+j}(v)\frac{%
H^{k-1+j}(v)}{h^{j-1}(v)}  \nonumber \\
& & -{r \choose 2}(k+r-1)_{(r-2)}\ _{r-1}M_k(v) \frac{h^{\prime k+1}(v)}{%
h^2(v)}+ \sum_{j=3}^{r} c_j(v)H^{k-1+j}(v)  \nonumber
\end{eqnarray}
where $c_j(v)$ are as in the statement of Lemma~4. For simplicity, further
on in the proof we drop the left subscript $r-1$ in $\ _{r-1}M_j(v)$ and
write $M_j(v)$ instead. Multiplying both sides of (\ref{first_eqn}) by $%
h^{r-1}(v)(v-l_F)^{r-1}/H^{k+r-1}(v)>0$ and making the change of variables
\[
w(v)=\frac{h(v)}{H(v)}(v-l_F),\qquad (v> l_F),
\]
we find (for simplicity we write $w$ for $w(v)$)
\begin{eqnarray}  \label{second_eqn}
(r-1)!g^{(k+r-1)}(v)w^{r-1}\frac{h(v)}{H(v)} & = & w^{r-1}\frac{h(v)}{H(v)}%
\sum_{j=0}^2 {r \choose j} (k+r-1)_{(r-j)}M_{k-1+j}(v)\frac{H^{j}(v)}{%
h^j(v)}  \nonumber \\
& & \hspace{-1cm}-{r \choose 2}(k+r-1)_{(r-2)}M_k(v)w^{r-2} \frac{%
h^{\prime }(v)}{h(v)}(v-l_F)+ S_1(v),
\end{eqnarray}
where
\[
S_1(v)=w^{r-1}\sum_{j=3}^{r} c_j(v)H^{j-1}(v).
\]
Referring to Lemma 1 with $n=r$ and $u=l_F$, we write (\ref{second_eqn}) as
\begin{eqnarray}  \label{third_eqn}
\lefteqn{\hspace{-1.5cm}w^{r-1}\frac{h(v)}{H(v)}\sum_{j=0}^2{r \choose
j}(k+r-1)_{(r-j)}(v-l_F)^jM_{k-1+j}(v)+ S_2(v)}  \nonumber \\
& = & w^{r-1}\frac{h(v)}{H(v)}\sum_{j=0}^2 {r \choose j}
(k+r-1)_{(r-j)}M_{k-1+j}(v)\frac{H^{j}(v)}{h^j(v)}  \nonumber \\
& & -{r \choose 2}(k+r-1)_{(r-2)}M_k(v)w^{r-2} \frac{h^{\prime }(v)}{h(v)}%
(v-l_F)+ S_1(v),
\end{eqnarray}
where
\[
S_2(v)=w^r\sum_{j=3}^{r}{r \choose j}%
(k+r-1)_{(r-j)}(v-l_F)^{j-1}M_{k-1+j}(v).
\]
It follows, from (\ref{third_eqn}), after simplifying and rearranging terms,
that
\begin{eqnarray*}
\lefteqn{\hspace{-2cm}w(w-1)(k+1)rM_k(v)+(w^2-1){r \choose
2}M_{k+1}(v)(v-l_F)} \\
& = & -{r \choose 2}\frac{h^{\prime }(v)}{h(v)}(v-l_F)M_k(v)+ \frac{%
S_1(v)-S_2(v)}{w^{r-2}(k+r-1)_{(r-2)}}.
\end{eqnarray*}
Finally, applying (\ref{frac_h}), we obtain
\begin{eqnarray}  \label{M_eqn}
\lefteqn{\hspace{-2cm}w(w-1)(k+1)rM_k(v)+(w^2-1){r \choose
2}M_{k+1}(v)(v-l_F)}  \nonumber \\
& = & -{r \choose 2}\left[\frac{w^{\prime }(v-l_F)}{w}+w-1\right]M_k(v) +
\frac{S_1(v)-S_2(v)}{w^{r-2}(k+r-1)_{(r-2)}}.
\end{eqnarray}

If $F$ is exponential, then $w(v)\equiv 1$. Since the exponential $F$ given
by (\ref{exp_type}) satisfies (\ref{exp_property}), we have that $w(v)\equiv
1$ is a solution of the above equation. To complete the proof we must show
that $w(v)\equiv 1$ is the only solution of (\ref{M_eqn}). Suppose $w(v)$ is
a solution of (\ref{M_eqn}) and there exists a value $v_1$ such that $%
w(v_1)\ne 1$ and $v_1>l_F$. We want to reach a contradiction. Since $F$ is
twice differentiable, we have that $w(v)$ is continuous with respect to $v$
and hence $w(v)\ne 1$ for $v$ in an open interval around $v_1$. (For a
similar argument see Lemma 3 in Su et al. (2008).) Let
\[
\label{a_inf} v_0=\inf\{v|w(v)\ne 1\}.
\]
Since, by (\ref{fraction_lim_1}), $w(l_F+)=1$, we have $v_0\ge l_F$. We
shall prove that $v_0=l_F$. Assume on contrary that $v_0>l_F$. Then $w(v)=1$
if $l_F<v\le v_0$ and integration of (\ref{new_var}) implies that $h(v)$ is
constant-valued in this interval. This contradicts the assumption (ii).
Therefore $v_0=l_F$ and hence equation (\ref{M_eqn}) holds for all $v>l_F$.
Dividing (\ref{M_eqn}) by $w-1\ne 0$, we obtain
\begin{eqnarray}  \label{M_eqn2}
\lefteqn{\hspace{-2cm}w(k+1)rM_k(v)+(w+1){r \choose 2}M_{k+1}(v)(v-l_F)}
\nonumber \\
& = & -{r \choose 2}\left[\frac{w^{\prime }(v-l_F)}{w(w-1)}+1\right]M_k(v)
+ \frac{S_1(v)-S_2(v)}{w^{r-2}(w-1)(k+r-1)_{(r-2)}}.
\end{eqnarray}
Passing to the limit as $v\to l_F+$ in the left-hand side of (\ref{M_eqn2}),
we find
\begin{eqnarray}  \label{lhs}
\lefteqn{\lim_{v\to l_F+}\left[wr(k+1)\ _{r-1}M_k(v)+(w+1){r \choose 2}\
_{r-1}M_{k+1}(v)(v-l_F)\right]}  \nonumber \\
& = & r(k+1)\lim_{v\to l_F+}w\ _{r-1}M_k(v)+{r \choose 2}\lim_{v\to
l_F+}(w+1)(v-l_F)\ _{r-1}M_{k+1}(v) \\
& = & r(k+1)\ _{r-1}M_{k}(l_F+),  \nonumber
\end{eqnarray}
where by (\ref{M_at_l_F_extra}), $\lim_{v\to l_F+}(v-l_F)\ _{r-1}M_{k+1}(v)=0
$ provided that $|g^{(k+r)}(l_F+)|<\infty$.

Now we turn to the right-hand side of (\ref{M_eqn2}). First, consider the
case $r=2$. Since $S_1(v)=S_2(v)=0$, we have for the right-hand side of (\ref%
{M_eqn2})
\begin{equation}  \label{rhs2}
\lim_{v\to l_F+}-\left[\frac{w^{\prime }(v-l_F)}{w(w-1)}+1\right]\
_1M_k(v)=- 2\ _1M_k(l_F+)
\end{equation}
where by the second equation in (\ref{fraction_lim_1}), $\lim_{v\to
l_F+}(v-l_F)w^{\prime }/(w-1)=1$. The equations (\ref{lhs}) and (\ref{rhs2})
imply $2(k+1)\ _1M_k(l_F+)=-2\ _1M_{k}(l_F+)$, which is not possible. This
proves that $w(v)\equiv 1$ is the only solution of  (\ref{M_eqn}) when $r=2$.

Let $r\ge 3$. Consider
\begin{equation}  \label{lim_1}
\lim_{v\to l_F+}\frac{S_2(v)}{w-1}=\lim_{v\to l_F+}w^r\frac{(v-l_F)^2}{w-1}%
\sum_{j=3}^{r}{r \choose j}(k+r-1)_{(r-j)}(v-l_F)^{j-3}\
_{r-1}M_{k-1+j}(v).
\end{equation}
By (\ref{m_frac}) we have $\lim_{v\to l_F+}(v-l_F)^2/(w-1)=0$. In addition,
by Lemma 2 we have that if $|g^{(k+r+2)}(l_F+)|<\infty$, then
\[
\limsup_{v\to l_F+}|\ _{r-1}M_{k+2}(v)|<\infty \quad \mbox{and}\quad
\lim_{v\to l_F+}(v-l_F)^{j-3}\ _{r-1}M_{k-1+j}(v)=0, \quad j=4, 5, \ldots r.
\]
Therefore, under the assumptions of the theorem, the limit in (\ref{lim_1})
is zero.

Let us now prove that 
\begin{equation}  \label{lim_2}
\lim_{v\to l_F+}\frac{S_1(v)}{w-1}=\lim_{v\to l_F+}\frac{H^2(v)}{w-1}%
\sum_{j=3}^{r} c_j(v)H^{j-3}(v)=0.
\end{equation}
It is not difficult to see that $\lim_{v\to \l _F+}H^{2}(v)/(w-1)=0$. Indeed,
Assumption (ii), (\ref{new_var}), and the first part of (\ref{fraction_lim_1}) together imply that $H(v)\sim const.(v-l_F)$, where $const.$ is not zero.
The limit assertion now follows from (\ref{m_frac}). Hence, to prove (\ref{lim_2}), it is sufficient to establish that the sum
in its right-hand side is finite. According to (\ref{d_formula2}) with $n=r$%
, this is true if $h(l_F+)\ne 0$, $|h^{(r-1)}(l_F+)|<\infty$ and $%
|g^{(k+2r-1)}(l_F+)|<\infty$, which hold by the assumptions of the theorem.

Taking into account (\ref{lim_1}) and (\ref{lim_2}), passing to the limit in
(\ref{M_eqn2}) as $v\to l_F$, we obtain, similarly to the case $r=2$, that $%
r(k+1)\ _{r-1}M_{k}(l_F+)=-r(r-1)\ _{r-1}M_{k}(l_F+)$ as $v\to l_F+$. This
contradiction proves that $w(v)\equiv 1$ is the only solution of  (\ref%
{M_eqn}). The proof of the theorem is complete.

\begin{acknowledgements}
I thank both referees for their very helpful critique and suggestions.
\end{acknowledgements}

\end{document}